\documentclass[12pt]{article}

\usepackage{amsfonts,latexsym,amssymb,hyperref}

\begin{document}

\newtheorem{definition}{\bf{Definition}}[section]
\newtheorem{lemma}{\bf{Lemma}}[section]
\newtheorem{proposition}{\bf{Proposition}}[section]
\newtheorem{theorem}{\bf{Theorem}}[section]
\newtheorem{remark}{\sc{Remark}}[section]

\title{A limit model for thermoelectric equations}
\author{Luisa Consiglieri\footnote{Independent Researcher Professor, Portugal.
\href{http://sites.google.com/site/luisaconsiglieri}{http://sites.google.com/site/luisaconsiglieri}
}}
\date{}
\maketitle

\begin{abstract}
We analyze the asymptotic behavior corresponding to the arbitrary
high conductivity of the heat in the thermoelectric devices.
This work deals with a steady-state multidimensional thermistor problem, considering 
  the Joule effect and both spatial and temperature dependent transport
coefficients under some real
boundary conditions in accordance with the Seebeck-Peltier-Thomson cross-effects.
Our first purpose is that the existence of a weak solution 
  holds true under minimal assumptions on the data, as 
in particular  nonsmooth domains.
Two existence results are studied under different assumptions on the electrical conductivity. 
 Their proofs are based on a fixed point argument,
compactness methods, and
existence and regularity  theory for elliptic scalar equations.
The second purpose is to show  the existence of a limit model illustrating  the
asymptotic situation. 
\end{abstract}



\section{Introduction}

In the recent years, the mathematical study of the thermistor problem has been of interest by various
authors (see for instance \cite{ac,cim,com,mg,xu} and the references therein) however the cross effects are neglected.
Here we study a mathematical model for thermoelectric conductors,
introduced in \cite{zamm}, taking into account
  the presence of the Seebeck-Peltier-Thomson and the Joule effects.
Indeed, in the thermodynamics analysis, the Joule effect is given by $|{\bf j}|^2/
\sigma(\cdot,\theta)$  with (cf. \cite{j})
\begin{equation}
{\bf j}=- 
\sigma(\cdot,\theta)(\alpha(\cdot,\theta)\nabla\theta+\nabla\phi),\label{dj}
\end{equation} 
representing the current density,
 $\theta$ denoting the temperature,
$\phi$ is the electric potential, 
and  the electrical conductivity $\sigma$ is a known positive function.
The electrical conductivity
is assumed temperature
 dependent and this is different at different places  along the material
due to the molecular structure.

The Seebeck coefficient $\alpha$ is a given nonlinear function, dependent both in space and temperature,  with constant sign
observing that
the sign of the Seebeck coefficient 
corresponds to the sign of the Hall effect \cite{ioffe}. 
The Thomson effect is $-\partial\alpha/\partial T(\cdot,\theta)\theta\nabla\theta\cdot{\bf j}$,
where $\partial/\partial T$ means the derivative with respect to the real variable.
Due to the first Kelvin relation,
$\pi(\theta)=\theta\alpha(\theta)$,
 $\theta\alpha(\theta)$ corresponds  to the Peltier coefficient.
Due to the second Kelvin relation, $\mu(\theta)=\theta\alpha'(\theta)$,
this coefficient is known as  Thomson coefficient. Although $\mu$ is the only thermoelectric coefficient
directly measurable for individual materials \cite{gasser,km},
and the Seebeck coefficient appears as $\alpha(T)=\int_{T_r}
^T {\mu(t)\over t}dt$ for some reference temperature $T_r$,
we keep the Seebeck coefficient as a given function as it is usual in the literature 
for the thermistor problem when this cross-effect coefficient is taken into account \cite{badii,cim2,xu2}.
Thus, the resulting PDE's system is strongly coupled.

The thermoelectric problem under study reads (see its derivation in \cite{zamm,j})

\noindent ($\mathcal P$)
 Find the pair
temperature-potential $(\theta,\phi)$ such that
\begin{eqnarray}
-\nabla\cdot(k\nabla\theta)=
\sigma(\cdot,\theta)\alpha(\cdot,\theta)(\alpha(\cdot,\theta)+{\partial\alpha\over\partial T}(\cdot,\theta)\theta)|\nabla\theta|^2
+\nonumber\\+
\sigma(\cdot,\theta)(2\alpha(\cdot,\theta)+{\partial\alpha\over\partial T}(\cdot,\theta)\theta)\nabla\theta\cdot
\nabla\phi+
\sigma(\cdot,\theta)|\nabla\phi|^2+g
&&\mbox{ in }\Omega,\label{pbu}\\
\label{defphi}
-\nabla\cdot(\sigma(\cdot,\theta)\nabla\phi)=\nabla\cdot
\left({ \sigma(\cdot,\theta)}\alpha(\cdot,\theta)\nabla\theta\right)
&&\mbox{ in }\Omega,
\\
k\nabla \theta\cdot{\bf n}=-\alpha(\cdot,\theta)\theta h
&&\mbox{ on }\partial\Omega,
\label{kh}\\
\sigma(\cdot,\theta)(\nabla\phi+\alpha(\cdot,\theta)\nabla\theta)\cdot{\bf n}=h
&&\mbox{ on }\partial\Omega,\label{pbn}
\end{eqnarray}
where $\Omega $ is a convex bounded domain of $\mathbb R^n$ $(n\geq 2)$,
 and it may represent  electrically conductive rigid solids
such as for instance  thermistors, thermocouples,
 resistive thermal devices (also called resistance temperature detectors) or thermoelectric coolers (see 
 \cite{bulu,chakra,yama} and the references therein).
 Notice that the boundary $\partial\Omega$ is  Lipschitz
since every bounded convex open subset of $\mathbb R^n$ always has a Lipschitz boundary \cite[Section 1.2]{grisv}.

Here $\bf n$ is the unit outward  normal to the boundary $\partial\Omega$,
 $g$ denotes the  external heat sources and $h$ denotes the surface current source.
The thermal conductivity $k$  is the known positive coefficient of 
the Fourier law.
Finally, we remark that (\ref{kh}) is known as  the linear Newton
law of cooling with $\alpha(\theta)h$ representing the heat transfer coefficient.

Our two main purposes are: (1)
to improve the existence result stated in \cite{zamm}
for three dimensional space, i.e.
the solution belongs to the Sobolev space $W^{1,p}(\Omega)$ with $p>n=3$.
Indeed, we state the existence result for the multidimensional case
since its proof is valid for any $n\geq 2$;
(2) to pass to the limit on the thermal
conductivity, $k\rightarrow+\infty$,  in
order to show the existence of the limit model
\begin{eqnarray}\label{t1}
\nabla\cdot (\sigma(\cdot,\Theta)\nabla\phi)=0 \quad\mbox{ in }\Omega;\\
\sigma(\cdot,\Theta)\nabla\phi\cdot{\bf n}=h \quad\mbox{ on }\partial\Omega,\label{t2}
\end{eqnarray}
for some positive constant $\Theta$  solving an implicit scalar equation.
The nonlocal problem (\ref{t1})-(\ref{t2}) is known as the shadow system,
the heuristic designation introduced by Nishiura \cite{n} in order
to exhibit minimal dynamics  displaying the mechanism of basic pattern formation.
Here the shadow system turns out 
the dynamic relation among the trivial rest states (constant solutions) and the large amplitude voltages.
The present result generalizes  the homogeneous Neumann
boundary value problem already studied for nonlocal elliptic
problems \cite{cr}.

The regularity assumptions about the domain could be weakened if different techniques are provided
(see \cite{bw,dk,dauge,elsch,med,zan} and the references therein).
We refer to \cite{costa} the existence of singularities of electromagnetic
fields at corners and edges of a bounded Lipschitz domain with piecewise plane boundary.

The contents of this work are as follow.
In next Section we state  the assumptions and the main results.
 Section \ref{aux} deals with existence and regularity results for auxiliary problems.
We prove in Section \ref{sexist} the  existence result for $n\geq 2$
when the electrical conductivity is assumed be a uniformly continuous
function, and
we obtain in Section \ref{sexist1} its limit model.
The proof of the existence result valid for $n= 2$ when the electrical conductivity is assumed be discontinuous
 on the space variable and the corresponding asymptotic limit model are postponed in Section \ref{2d}.

\section{Assumptions and main results}

In order to establish the existence results we assume the following set of hypotheses on the data.
\begin{description}
\item[(H1)] $k$ is a positive constant.

\item[(H2)] $\sigma:\Omega\times\mathbb R\rightarrow\mathbb R$ is a
  Carath\'eodory function, i.e. measurable with respect to $x\in\Omega$ and
  continuous with respect to $T\in\mathbb R$, and furthermore
\begin{equation}
\exists \sigma_\#,\sigma^\#>0:\quad
\sigma_\#\leq \sigma(\cdot,T)\leq \sigma^\#,\quad\forall T\in \mathbb R,\
\mbox{a.e. in }\Omega.\label{smin}\end{equation}

\item[(H3)] $ g$ belongs to the Lebesgue space $ L^{p/2}(\Omega)$, $p>2$, and
 $h\in C(\partial\Omega)$ is such that verifies the compatibility condition
\begin{equation}\label{cc}
\int_{\partial\Omega}h\ ds=0,
\end{equation}
where $ds$ represents the element of  surface area.

\item[(H4)] $\alpha\in C^{}(\bar\Omega\times\mathbb R)$ is such that
 \begin{equation}
\exists L_\alpha>0:\quad|{\alpha}(x_1,T_1)-{\alpha}(x_2,T_2)|\leq L_\alpha(|x_1-x_2|+|T_1-T_2|),\quad
\label{la}
\end{equation}
for all $x_1,x_2\in\bar\Omega$ and for all $T_1,T_2\in \mathbb R$, with $|x|$ representing the euclidean
norm and $|T|$ the absolute value of a real number.
Moreover, for all $x\in\Omega$,
the mapping $T\mapsto\alpha(x,T)$ is continuously differentiable in $\mathbb R$ and its derivative satisfies
\begin{equation}
\exists\mu^\#>0:\quad|{\partial\alpha\over\partial T}(x,T)|\leq
\left\{\begin{array}{ll}
\mu^\#,\quad&\mbox{ if }|T|\leq 1\\
\mu^\#/|T|,\quad&\mbox{ if }|T|>1 .\end{array}\right.
 \label{alphat}
 \end{equation}
The following two different cases will be addressed:
\begin{description}
\item[(H4)$_+$ for materials with positive $\alpha$:]
\begin{equation}
\exists\alpha_\#,\alpha^\#>0:\quad
\alpha_\#\leq \alpha(x,T)\leq \alpha^\#,\quad\forall x\in \bar\Omega,\quad\forall T\in \mathbb R.\label{amm}
\end{equation}
Moreover, there exist an open subset $\Gamma\subset \partial\Omega$ such that meas$(\Gamma)>0$
 and meas$(\partial \Omega\setminus\bar\Gamma)>0$
and constants $h_\#>0$ and $h_1<0$ such that
\begin{eqnarray}\label{defh}\qquad
h_1>-{\min\{k,\alpha_\#h_\#\}\over C_1\alpha^\#},
\quad h(x)\geq
\left\{\begin{array}{ll}
 h_\#, \quad&\mbox{a.e. } x\in \Gamma\\
h_1, &\mbox{a.e. }x\in\Sigma,
\end{array}\right.
\end{eqnarray}
where
 $\Sigma:=\partial \Omega\setminus\bar\Gamma$ and
 $C_1$ denotes the continuity constant of the embedding $H^1(\Omega)\hookrightarrow L^2(\Sigma)$, i.e.
 \begin{equation}\label{cc1}
\|\theta\|_{2,\Sigma}^2\leq C_1\left(\|\nabla\theta\|_{2,\Omega}^2+\|\theta\|_{2,\Gamma}^2\right)
\end{equation}
for every $\theta\in H^1(\Omega)$;
 
\item[(H4)$_-$ for materials with negative $\alpha$:]
\begin{equation}
\exists\alpha_\#,\alpha^\#>0:\quad
- \alpha^\#\leq\alpha(x,T)\leq -\alpha_\#,\quad\forall x\in \bar\Omega,\quad\forall T\in \mathbb R.\label{amm-}
\end{equation}
Moreover, there exist an open subset $\Gamma\subset \partial\Omega$ such that meas$(\Gamma)>0$
 and meas$(\partial \Omega\setminus\bar\Gamma)>0$
and constants $h_\#>0$ and $h_1<0$ such that
\begin{eqnarray}\label{defh-}\qquad
h_1>-{\min\{k,\alpha_\#h_\#\}\over C_1\alpha^\#},
\quad h(x)\leq
\left\{\begin{array}{ll}
- h_\#, \quad&\mbox{a.e. } x\in \Gamma\\
-h_1, &\mbox{a.e. }x\in\Sigma.
\end{array}\right.
\end{eqnarray}
\end{description}

\end{description}

\begin{remark}
In particular, (\ref{cc1}) holds for the unity function then we have
\begin{equation}\label{c1}
|\Sigma|\leq C_1|\Gamma|,
\end{equation}
where $|\Gamma|=$meas$(\Gamma)$ and $|\Sigma|$=meas$(\partial \Omega\setminus\bar\Gamma)$.
\end{remark}

{\bf Some remarks on the assumptions.}
\begin{enumerate}
\item The heat conductivity can be an uniformly continuous function on both variables verifying
\[\exists k_\#,k^\#>0:\quad
k_\#\leq k(x,T)\leq k^\#,\quad\forall x\in \Omega,\quad\forall T\in \mathbb R.
\]
Since our purpose is to study the asymptotic behavior as $k\rightarrow \infty$,
we assume it as constant.

\item The case $\Gamma=\partial\Omega$ is excluded, from the fact that the Gauss theorem yields
the necessary condition (\ref{cc}) of the existence of a solenoidal function satisfying (\ref{dj}) and (\ref{pbn})
 (cf. (\ref{defphi})).

\item 
From the assumption (\ref{alphat}) we obtain that
$|{\partial\alpha\over\partial T}(x,T)|\leq \mu^\#,$
for all $x\in \Omega$ and for all $T\in\mathbb R$,
and that the Thomson coefficient is bounded, i.e.
$|\mu(x,T)|=|T{\partial\alpha\over\partial T}(x,T)|\leq \mu^\#$, for all $x\in \Omega$ and for all $T\in \mathbb R$.
\item By the Weierstrass Theorem, any continuous function defined on a compact set (of
$\mathbb R^{n+1}$) is bounded. Then the upper bound in
(\ref{amm}) as well as the lower bound in (\ref{amm-}) could be given reduced to $|T|>1$.
\end{enumerate}

We define the Nemytskii operators
\begin{eqnarray*}
\sigma(\theta)=\sigma(\cdot,\theta(\cdot));\qquad
\alpha(\theta)=\alpha(\cdot,\theta(\cdot));\qquad
{\partial\alpha\over\partial T}(\theta)={\partial\alpha\over\partial T}(\cdot,\theta(\cdot)),
\end{eqnarray*}
that map $L^1(\Omega)$ into $L^q(\Omega)$, for all $q<\infty$.
Their designation is kept in order to clarify the presentation.

We endow the Sobolev space $W^{1,p}(\Omega)$, $p>1$, with the equivalent norms:
\begin{itemize}
\item for the temperature solution
\[
\|\theta\|_{1,p,\Omega}=\|\nabla\theta\|_{p,\Omega}+\|\theta\|_{p,\Gamma};\]
\item for the potential solution
\[\|\phi\|_{1,p,\Omega}=\|\nabla\phi\|_{p,\Omega},\]
considering the correspondent Poincar\'e inequality.
\end{itemize}

\begin{definition}\label{d1}
We say that $(\theta,\phi)$ is a {\bf weak solution} to
(\ref{pbu})-(\ref{pbn}) if $(\theta,\phi)\in W^{1,p}(\Omega)^2$, for $p> n$,
 and it satisfies 
\begin{eqnarray}
k\int_\Omega \nabla\theta\cdot\nabla\eta dx
+\int_{\partial\Omega}\alpha(\theta)h\theta \eta ds=\nonumber\\
=
\int_\Omega \sigma(\theta)\Big(\alpha(\theta)(\alpha(\theta)+
{\partial\alpha\over\partial T}(\theta)\theta)|\nabla\theta|^2+(2\alpha(\theta)+{\partial\alpha\over\partial T}(\theta)\theta)\nabla\theta\cdot\nabla\phi\Big)\eta dx
\nonumber\\
+\int_\Omega \left({\sigma(\theta)}|\nabla\phi
|^2+g\right)
\eta dx,\quad\forall \eta \in W^{1,p'}(\Omega);\qquad\label{heatj}\\
\int_\Omega\sigma(\theta)\nabla\phi\cdot \nabla\eta dx=\nonumber\\
=-\int_\Omega\sigma
(\theta)\alpha(\theta)\nabla\theta\cdot\nabla\eta dx
 +\int_{\partial\Omega}
h\eta ds,\quad\forall \eta\in W^{1,p'}(\Omega),\qquad\label{pbphij}
 \end{eqnarray}
where $p'$ denotes the conjugate exponent to $p$, $p'=p/(p-1)$.
\end{definition}

\begin{remark}
\label{mpnl}
If $\Omega$ is of class $C^{0,1}$ and
$mp>n,$ then  the Morrey-Sobolev 
embedding holds
\[W^{m,p}(\Omega)\hookrightarrow C^{0,m-n/p}(\bar\Omega).\]
The vector field $\bf n$  belongs only to ${\bf L}^\infty(\partial\Omega)$ if it
is the unit outward normal vector
to the boundary of $C^{0,1}$ domains.
When $p>{n},$ the embedding $W^{1,p'}
(\Omega)\hookrightarrow L^{pn/(pn-n-p)}(\Omega)$
is valid. Thus the $L^{p/2}$ behavior
of the quadratic terms $|\nabla\theta|^2$, $\nabla\theta\cdot\nabla\phi$,
$|\nabla\phi|^2$ is meaningful
on the right hand side of (\ref{pbu}) since $p/2>pn/(p+n)$.\end{remark} 

Let us extend the existence results whose can be found in \cite{zamm}.
The first main theorem states the existence of weak solutions to the problem under study, 
strengthening the assumption (H2), i.e. strengthening the regularity on $\sigma$.
\begin{theorem}\label{exist}
Assume $n\geq 2$, (H1) and (H3)-(H4) hold.
Additionally we assume that $
\sigma:\Omega\times\mathbb R\rightarrow\mathbb R$ is a uniformly continuous function satisfying (\ref{smin}) and the smallness condition (\ref{ss}) is satisfied. 
 Then the variational problem (\ref{heatj})-(\ref{pbphij}) admits a weak
solution $(\theta,\phi)\in W^{2,pn/(p+n)}(\Omega)\times W^{1,p}(\Omega)$,
with $p>n$, in the sense of Definition \ref{d1}, such that $\int_{\partial\Omega}\phi ds$.
Moreover, the following estimates hold
\begin{eqnarray}
\|\nabla\phi\|_{2,\Omega} 
&\leq& {1\over\sigma_\#}(\sigma^\#\alpha^\#\|\nabla\theta\|_{2,\Omega}+C_2\|h\|_{2,\partial\Omega} ) ;
\label{cotaphi}\\
\label{cotaup}
\|\theta\|_{1,p,\Omega}&\leq &C_p;\\
k\|\nabla \theta\|_{2,\Omega}^2&\leq& C_3\Big(\alpha^\#|h_1|C_p^2+\nonumber\\
\label{cotau}
&&+(AC_p^2+BC_p
\|\nabla\phi\|_{2,\Omega}+\sigma^\#\|\nabla\phi\|_{2,\Omega}^2+\|g\|_{p/2,\Omega})C_p\Big),
\end{eqnarray}
where
 $C_2$ denotes the Poincar\'e-Sobolev continuity constant of the embedding $H^1(\Omega)\hookrightarrow L^2(\partial\Omega)$,
  $C_3$ denotes the continuity constant of the embedding $W^{1,p}(\Omega)\hookrightarrow C(\bar\Omega)$,
 \begin{equation}\label{ab}
A:=\sigma^\#\alpha^\#(\alpha^\#+\mu^\#),\qquad
B:=\sigma^\#(2\alpha^\#+\mu^\#);
\end{equation}
 and $C_p$ denotes a positive constant independent on $k$ if $k>\alpha_\#h_\#$.
\end{theorem}
 If we assume $\sigma\in C(\bar\Omega\times\mathbb R)$ as in \cite{zamm}, 
  the Nemytskii operator maps $C(\bar\Omega)$ into $C(\bar\Omega)$ which implies the uniform
continuity on the spatial variable. 

Next, we establish the existence of a solution  $\phi=\phi(\Theta)$ to
(\ref{t1})-(\ref{t2}) where $\Theta$ is solution of an implicit scalar equation.
\begin{theorem}\label{exist1}
Under the assumptions of Theorem \ref{exist}, there exist weak
solutions $(\theta_k,\phi_k)\in W^{2,pn/(p+n)}(\Omega)\times W^{1,p}(\Omega)$,
with $p>n$, to the variational problem (\ref{heatj})-(\ref{pbphij}),
such that
\begin{eqnarray*}
\theta_k\rightarrow  \Theta \quad\mbox{ in }
H^1(\Omega), \qquad \phi_k\rightharpoonup\phi \mbox{ in $H^1(\Omega)$,}
\end{eqnarray*}
with $\Theta$  solving the implicit scalar equation
\begin{equation}\label{tt}
\Theta\int_{\partial\Omega}\alpha(\cdot,\Theta)hds
=\int_\Omega\sigma(\cdot,\Theta)|\nabla\phi|^2dx+\int_{\Omega}gdx
\end{equation}
and $\phi$ solving (\ref{t1})-(\ref{t2}).
Moreover 
\begin{equation}\label{tt1}
{\sigma_\#\int_\Omega
|\nabla\phi|^2dx+\int_{\Omega}g dx \over\alpha^\#\int_{\partial\Omega}|h|ds} \leq
\Theta\leq{\sigma^\#\int_\Omega
|\nabla\phi|^2dx+\int_{\Omega}gdx \over \alpha_\#h_\#|\Gamma|+\alpha^\#h_1|\Sigma|} .
\end{equation}
\end{theorem}

\begin{remark}
In face of (\ref{cc}) if the Seebeck coefficient is only a function on the temperature then (\ref{tt}) reads
$$
0=\int_\Omega\sigma(\cdot,\Theta)|\nabla\phi|^2dx+\int_{\Omega}gdx.$$
The presence of a generic heat source
invalids the limit model (\ref{t1})-(\ref{t2}) connected with the original thermoelectric problem
introduced in \cite{zamm}.
\end{remark}

The above results can be proved if the convexity of $\Omega$ is replaced by weaker assumptions, 
for instance
 when $\Omega$ is a plane bounded domain
with Lipschitz and piecewise $C^2$ boundary whose angles are all convex 
\cite[p. 151]{grisv},
or when $\Omega$ is a plane bounded domain
with curvilinear polygonal $C^{1,1}$ boundary whose angles are all strictly convex \cite[p. 174]{grisv}.
For general bounded domains $\Omega$ of $\mathbb R^n$ with Lipschitz  boundary, it is known that
the integrability exponents for the gradients
 of the potential and temperature solutions may be larger than 3 \cite{elsch,zan},
if the restriction to the case of
uniformly continuous coefficients in (\ref{heatj})-(\ref{pbphij}) is assured. However,
 a generalization for such nonsmooth domains of Theorem \ref{exist}
and its limit model is not a direct consequence.
Indeed, new proofs will be needed because the compact
embedding 
$W^{2,pn/(p+n)}(\Omega)\hookrightarrow\hookrightarrow W^{1,p}(\Omega)$
is crucial to provide the weak continuity of the operator in the fixed point argument.

For the two-dimensional limit model, let us show the existence result
under the minimal regularity on $\sigma$.
\begin{theorem}\label{teo1}
Assume $n= 2$ and (H1)-(H4) hold. Then the problem 
(\ref{pbu})-(\ref{pbn}) has a weak solution in the sense of Definition \ref{d1},
 for some $p> 2=n$, under sufficiently small data.
\end{theorem}

\begin{proposition}\label{exist2}
Under the assumptions of Theorem \ref{teo1}, there exists
$\phi$ solving (\ref{t1})-(\ref{t2}) with $\Theta$  solving the implicit scalar equation (\ref{tt}).
\end{proposition}

The study of the existence of three-dimensional weak solutions to the variational problem  (\ref{heatj})-(\ref{pbphij}),
under the assumption that the mapping 
$x\in\Omega\mapsto \sigma(x,T)$ is discontinuous for every $T\in \mathbb{R}$, is still an open problem.
When the coefficient $\sigma$ of the principal part of the divergence form elliptic
equation (\ref{defphi}) is a discontinuous function  on the spatial variable,
it invalidates the smoothness of the solution as is carried out in the literature \cite{cfl,kp,lsw,mey,rag}. 

\section{Auxiliary results}
\label{aux}

The existence of a unique solution $\phi\in H^1(\Omega)$, such that
$\int_{\partial\Omega} \phi ds=0$, to an auxiliary problem is consequence of 
Lax-Milgram lemma, for details see 
 \begin{proposition}[{\cite[Theorem 4.1]{zamm}}]\label{twp2}
 Let the assumptions (H2) and (\ref{amm}) 
 or (\ref{amm-}) be fulfilled.
 Assume  that $n\geq 2$, $\Omega\in C^{0,1}$,
 $\xi\in H^1(\Omega)$ and  $h\in  L^{p}(\partial\Omega)$  verify (\ref{cc}) for $p> 2(n-1)/n$.
Then there exists a unique weak solution $\phi\in H^1(\Omega)$, such that  $\int_{\partial\Omega}\phi ds=0$,
 to the variational problem
 \begin{equation}\label{auxphi}
\int_\Omega\sigma(\xi)\nabla\phi\cdot \nabla\eta dx=-\int_\Omega\sigma(\xi)\alpha(\xi)\nabla\xi\cdot\nabla\eta dx
 +\int_{\partial\Omega} h\eta ds,
\end{equation}
 for all $\eta \in H^1(\Omega)$ and in particular for all $\eta\in W^{1,p'}(\Omega)$.
\end{proposition}

Next, we establish some regularity for the potential auxiliary solution.
\begin{proposition}\label{twp1}
Let  $p>n$, $h\in L^p(\partial\Omega)$ verify (\ref{cc}),  
$\xi\in W^{2,q}(\Omega)$ with $q=pn/(p+n)$, (H4) be fulfilled, and $\phi\in H^{1}(\Omega)$ solve the problem 
 (\ref{auxphi}). If  $
\sigma:\Omega\times\mathbb R\rightarrow\mathbb R$ is a uniformly continuous function
verifying (\ref{smin}), then $\nabla\phi\in {\bf L}^{p}(\Omega)$ and it verifies
\begin{equation}
\|\nabla\phi\|_{p,\Omega} 
\leq C(n,p,\Omega,\sigma_\#,\sigma^\#)(\sigma^\#\alpha^\#\|\nabla\xi\|_{p,\Omega}+\|h\|_{p,\partial\Omega} ) 
.\label{cotap}
\end{equation}
\end{proposition}
{\sc Proof.}
For $p>n$, we have $q>n/2$,
\begin{eqnarray*}
\xi\in W^{2,q}(\Omega)\hookrightarrow C^{0,2-n/q}(\bar\Omega),\quad
\alpha(\xi) \nabla\xi\in {\bf W}^{1,q}(\Omega)
\hookrightarrow{\bf L}^{q(n-1)/(n-q)}(\partial\Omega),\\
\alpha(\xi)\nabla\xi\cdot{\bf n}\in { L}^{q(n-1)/(n-q)}(\partial\Omega)
\equiv  { L}^{p(n-1)/n}(\partial\Omega)\hookrightarrow {W}^{-1/p,p}(\partial\Omega), 
\end{eqnarray*}
for ${\bf n}\in{\bf L}^\infty(\partial\Omega)$.
Moreover,  $h/\sigma(\xi)\in W^{-1/p,p}(\partial\Omega)$ and
$\nabla\xi\in {\bf L}^p(\Omega)$ implies that $\nabla\cdot(\sigma(
\xi)\alpha(\xi)\nabla\xi)\in
\left({\bf W}^{1,p'}(\Omega)\right)'$ for $p>n$.
By appealing to the regularity theory \cite{daut1,grisv} for the solution $\phi\in H^{1}(\Omega)$
 of the boundary value problem (in the sense of distributions)
\begin{eqnarray*}
&&\nabla\cdot\left(\sigma(\xi)\nabla \phi+{ \sigma(\xi)}\alpha(\xi)\nabla\xi\right)=0\quad
\mbox{ in }\Omega\\
&&\sigma(\xi)(\nabla \phi+\alpha(\xi)\nabla\xi)\cdot{\bf n}=h\quad\mbox{on }\partial\Omega,
\end{eqnarray*}
and observing that $\xi\in W^{1,p}(\Omega)\hookrightarrow C^{0,1-n/p}(\bar\Omega)$ for $p>n$
warrants that $\sigma(\cdot,\xi)$ is  uniformly continuous, then
the  regularity of  weak solutions relative to $W^{1,p}(\Omega)$ and the estimate (\ref{cotap}) arise. $\Box$

The following result deals with the existence and uniqueness of a strong temperature  auxiliary solution.
 \begin{proposition}\label{propt}
Let  $p>n$,   $\xi\in W^{1,p}(\Omega)$, (H1) and (H3)-(H4) be fulfilled, and $\phi\in H^{1}(\Omega)$ solve the problem 
 (\ref{auxphi}). If  $\sigma:\Omega\times\mathbb R\rightarrow\mathbb R$ is a uniformly continuous function
verifying (\ref{smin}), then there exists a unique weak solution  $\theta\in
W^{2,pn/(p+n)}(\Omega)$  solving the problem, for all $\eta \in W^{1,p'}(\Omega)$,
\begin{eqnarray}
k\int_\Omega \nabla\theta\cdot\nabla\eta dx
+\int_{\partial\Omega}\alpha(\xi)h\theta \eta ds
=
\int_\Omega \left(F(\cdot,\xi,\nabla\xi,\nabla\phi)
+g\right)
\eta dx,
\label{auxu}\end{eqnarray}
with $F:\Omega\times\mathbb R^{2n+1}\rightarrow\mathbb R$ defined as 
$F(x,T,{\bf a},{\bf b})=$
\[\sigma(x,T)\Big(\alpha(x,T)(\alpha(x,T)+
{\partial\alpha\over\partial T}(x,T)T)|{\bf a}|^2+(2\alpha(x,T)+{\partial\alpha\over\partial T}(x,T)T){\bf a}
\cdot{\bf b}+|{\bf b}|^2\Big).\]
\end{proposition}
{\sc Proof.}
The existence and uniqueness of $\phi\in W^{1,p}(\Omega)$ is consequence of Propositions \ref{twp2} and \ref{twp1}.
By appealing to the elliptic equations theory \cite{grisv},
from $F(\xi,\nabla\xi,\nabla\phi)+g\in L^ {p/2}(\Omega)$,
the regularity theory for the Laplace equation in convex domains  guarantees the
 existence of a unique solution $\theta\in W^{2,pn/(p+n)}(\Omega)$ of    the Robin problem 
\begin{eqnarray*}
&&-k\Delta \theta=F(\cdot,\xi,\nabla\xi,\nabla\phi)+g\quad \mbox{ in }\Omega;\\
&&k\nabla \theta\cdot{\bf n}+\alpha(\xi) h\theta=0\quad \mbox{ on }\partial\Omega,
\end{eqnarray*}
taking into account that
the Korn perturbation method \cite[pp. 107-109]{grisv} can be adapted if
the coefficient   
 $\alpha(\cdot,\xi)h\in C(\partial \Omega)$ is such that  the assumption (\ref{la}) holds.  For this, we observe that
$\xi\in W^{1,p}(\Omega)\hookrightarrow C(\bar\Omega)$ and we recall (H3)-(H4). $\Box$

For the regularity of the potential  auxiliary solution $\phi$ when
it is the unique weak solution for Neumann problem to an elliptic second order equation in divergence form with
bounded and measurable coefficient, we can prove the following result.
\begin{proposition}\label{pphi}
If the assumptions of Proposition \ref{twp2} are fulfilled
with $ p=2$ and $\Omega$ is convex, then there
 exists a constant $\epsilon>0$ such that the weak solution $\phi\in H^1(\Omega)$
of (\ref{auxphi}) belongs to $W^{1,2+\epsilon}(\Omega)$, i.e.
\begin{equation}\label{cotaphie}
\|\nabla {\phi}\|_{2+\epsilon,\Omega}\leq K_2(\sigma^\#\alpha^\#\|\nabla\xi\|_{2,\Omega}+\|h\|_{2,\partial\Omega} ) ,
\end{equation}
with a constant $K_2>0$ only dependent on the data.
\end{proposition}
{\bf Proof.}
Denote the operator $A$ by
\[
\langle A\phi,\eta\rangle=\int_\Omega \sigma(\xi)\nabla\phi\cdot\nabla \eta dx.
\]
Then $\phi\in H^{1}(\Omega)$ is a weak solution to the second order elliptic differential equation
$Au=F$, under
\[F=\nabla\cdot( \sigma(\xi)\alpha(\xi)\nabla\xi)+h\in ( H^{1}(\Omega))'\hookrightarrow 
( W^{1,p}(\Omega))',\qquad \forall p\geq 2.\]
Since
the boundedness property 
\[
\sigma_\#\leq \sigma(\cdot,\xi)\leq \sigma^\#,
\mbox{ a.e. in }\Omega,\]
 is fulfilled,
considering that  $\xi\in L^1(\Omega)$ and the assumption  (\ref{smin}) on $\sigma$ holds,
then the Neumann  version of the general
 result on the  higher regularity for weak solutions to the  mixed boundary value
problems (cf. \cite[Theorem 1]{grog}, also \cite{grogr}) guarantees that
$\phi\in W^{1,2+\epsilon}(\Omega)$ for some $\epsilon >0$. $\Box$

Although Proposition \ref{pphi} is valid  for any dimensional space $(n\geq 2)$,
we only used it for $n=2$.  Let us precise its application in the following proposition.
 \begin{proposition}\label{propt2}
Let     $\xi\in W^{1,2+\epsilon}(\Omega)$, (H1)-(H4) be fulfilled, 
and $\phi\in H^{1}(\Omega)$ solve the problem 
 (\ref{auxphi}). 
If the assumptions of Proposition \ref{pphi} hold,
 then there exists a unique weak solution  $\theta\in
W^{2,2p/(p+2)}(\Omega)$  solving (\ref{auxu}) with $p=2+\epsilon$.
\end{proposition}
{\sc Proof.}
The imperative requirement of the embedding
$ W^{1,2+\epsilon}(\Omega)\hookrightarrow C(\bar\Omega)$ yields that
 $\alpha(\cdot,\xi)h\in C(\partial \Omega)$ provided by (H3)-(H4). Thus,
 Proposition \ref{pphi} ensures that the argument of the proof of Proposition \ref{propt}
 is still valid, concluding the claim.
$\Box$

\section{Proof of Theorem \ref{exist}}
\label{sexist}

First we recall the
Tychonoff extension to weak topologies of the Schauder fixed point
theorem \cite[pp. 453-456 and 470]{dsch}.
\begin{theorem}\label{fpt} 
Let $K$ be a nonempty   compact convex subset of a locally convex  space $X$. Let
${\mathcal L}:K\rightarrow K$ be a continuous operator. Then $\mathcal L$ has at least one fixed point.
\end{theorem}

If we provide any Banach space with the weak topology, every closed ball  is convex and weakly sequential compact.

In order to apply Theorem \ref{fpt},
let us consider the operator $\mathcal L$ defined in a closed ball
$\bar B_R \subset W^{2,pn/(p+n)}(\Omega)$ such that  
\[
\mathcal{L}:\xi\in\bar B_R\mapsto\phi\mapsto \theta\in
W^{2,pn/(p+n)}(\Omega),
\]
 where $\phi\in W^{1,p}(\Omega)$ solves the problem (\ref{auxphi}), for all $\eta \in W^{1,p'}(\Omega)$,
and $\theta$  solves the problem (\ref{auxu}).

The existence of a unique solution $\phi\in W^{1,p }(\Omega)$, such that
$\int_{\partial\Omega} \phi ds=0$, to the problem (\ref{auxphi}) is consequence of Propositions
\ref{twp2} and \ref{twp1}, and it  verifies (\ref{cotap}).

Hence,
for $p>n$, we find $\theta\in W^{1,p}(\Omega)$ from Proposition \ref{propt}, and the estimate
\begin{equation}
\| \theta\|_{2,pn/(p+n),\Omega}
\leq K(A\|\nabla\xi\|_{p,\Omega}^2+B\|\nabla\xi\|_{p,\Omega}
\|\nabla\phi\|_{p,\Omega}+\sigma^\#\|\nabla\phi\|_{p,\Omega}^2+
\|g\|_{p/2,\Omega}),
\label{cotatp}
\end{equation}
is verified with $K$ denoting a constant dependent on $\Omega$, $n$ and $p$,
 $A$ and $B$ given by (\ref{ab}), and $\varkappa:=\min\{k,\alpha_\#h_\#\}+\alpha^\# h_1C_1>0$.

 Thus $\mathcal L$ is well defined.
 
  Next, let us prove that
$\mathcal L(\bar B_R) \subset \bar B_R.$ Let $\xi\in \bar B_R$ be arbitrary
and $(\phi,\theta)$ be the corresponding solution solving
(\ref{auxphi}) and (\ref{auxu}). Thus (\ref{cotap}) and (\ref{cotatp}) read
\begin{eqnarray}
\|\phi\|_{1,p,\Omega} 
\leq {C}(\sigma^\#\alpha^\#R+\|h\|_{p,\partial\Omega} ) 
;\label{cotapr}\\
\label{cotatpr}
\| \theta\|_{2,pn/(p+n),\Omega}
\leq K(AR^2+BR
\|\nabla\phi\|_{p,\Omega}+\sigma^\#\|\nabla\phi\|_{p,\Omega}^2+
\|g\|_{p/2,\Omega}),
\end{eqnarray}
with $C=C(n,p,\Omega,\sigma_\#,\sigma^\#)$.
Inserting (\ref{cotapr}) into (\ref{cotatpr}) it follows
\[\| \theta\|_{2,pn/(p+n),\Omega}\leq a_2R^2+a_1R+a_0,
\]
where
\begin{eqnarray*}
a_2&=& K\sigma^\#\alpha^\#\left(1+{C\sigma^\#}\right)
\left(\alpha^\#\left(2+{C\sigma^\#}\right)+\mu^\#\right);
\\a_1&=& {KC\sigma^\#}\left(2\alpha^\#(1+ \sigma^\#)+\mu^\#\right)
\|h\|_{p,\partial\Omega};
\\a_0&=& K\left({C^2\sigma^\#}\|h\|_{p,\partial\Omega}^2+
\|g\|_{p/2,\Omega}\right).
\end{eqnarray*}
Therefore, $\mathcal L(\xi)=\theta\in\bar B_R$ if and only if $a_2R^2+(a_1-1)R+a_0\leq 0$,
i.e. for instance if the smallness condition 
\begin{equation}\label{ss}
a_1<1\quad\mbox{and}\quad 4a_0a_2<(1-a_1)^2
\end{equation}
is assumed. 

\subsection{The weak sequential continuity of $\mathcal L$}
\label{swc}

Let $\{\xi_m\}_{m\in\mathbb N}$ be a sequence in $\bar B_R$
 verifying
\begin{equation}\label{xim2}
\xi_m\rightharpoonup \xi\quad\mbox{in }
W^{2,pn/(p+n)}(\Omega)\hookrightarrow\hookrightarrow W^{1,p}(\Omega),
\end{equation}
and $({\phi}_m,\theta_m)$ is the correspondent solution to 
(\ref{auxphi}) and (\ref{auxu}),
 for each $m\in \mathbb N$. 
From the estimates (\ref{cotap}) and (\ref{cotatp}) we can extract a subsequence, still labeled
by $(\phi_m,\theta_m)$, such that
\[
\phi_m \rightharpoonup \phi\quad
\mbox{ in }W^{1,p}(\Omega),\qquad
{\theta}_m \rightharpoonup{\theta}\quad 
\mbox{ in }
W^{2,pn/(p+n)}(\Omega)\hookrightarrow\hookrightarrow W^{1,p}(\Omega).
\]
Thanks to Remark \ref{mpnl} it follows
\begin{equation}\label{xim1}
\xi_m \rightarrow \xi, \quad
{\theta}_m \rightarrow{\theta},\quad 
\phi_m\rightarrow\phi\quad \mbox{ in }C^{0,1-n/p}(\bar\Omega).
\end{equation}
In particular, $\int_{\partial\Omega}\phi_m ds=0\rightarrow \int_{\partial\Omega} \phi ds=0$.
By the continuity of the Nemytskii operators $\alpha$ and $\sigma$,
we can pass to the limit in (\ref{auxphi})$_m$ as $m$ tends to infinity,
 concluding that $\phi\in W^{1,p}(\Omega)$ is the limit solution, i.e. it verifies (\ref{auxphi}).

In the sequel, 
let us pass to the limit in (\ref{auxu})$_m$ as $m$ tends to infinity.
First,  the mapping $\xi\in L^1(\Omega)\mapsto
\alpha(\xi)\in L^r(\Omega),$ for all $r <+\infty,$ is continuous by
(H4), thus the passage to the limit to the left hand side of (\ref{auxu}) 
is straightforward. In order to study the RHS, we define
\begin{eqnarray*}
I_{1,m}&=&
\int_\Omega \sigma(\xi_m)\alpha(\xi_m)(\alpha(\xi_m)+
{\partial\alpha\over\partial T}(\xi_m)\xi_m)|\nabla\xi_m|^2\eta dx;\\
I_{2,m}&=&\int_\Omega\sigma(\xi_m)(2\alpha(\xi_m)+{\partial\alpha\over\partial T}(\xi_m)\xi_m)
\nabla\xi_m\cdot\nabla\phi_m\eta dx;\\
I_{3,m}&=&
\int_\Omega {\sigma(\xi_m)}|\nabla\phi_m|^2\eta dx.
\end{eqnarray*}
Recalling Remark \ref{mpnl}
we get $\eta\in W^{1,p'}
(\Omega)\hookrightarrow L^{pn/(pn-n-p)}(\Omega)\hookrightarrow L^{p/(p-2)}(\Omega)$ for $p>n$.

From (\ref{xim2}) we have
$|\nabla{\xi}_m|^2 \rightarrow |\nabla\xi|^2$ in ${ L}^{p/2}(\Omega)$.
Considering that  the mapping $\xi\in L^1(\Omega)\mapsto
\sigma(\xi)\alpha^2(\xi)\eta\in L^{p/(p-2) }(\Omega)$ is continuous
thus the first term in $I_{1,m}$ passes to the limit as $m$ tends to infinity.
Using (\ref{xim2}) and  (\ref{xim1}) we have
$\xi_m|\nabla{\xi}_m|^2 \rightarrow \xi|\nabla\xi|^2$ in ${ L}^{p/2}(\Omega)$.
Considering that  the mapping $\xi\in L^1(\Omega)\mapsto
\sigma(\xi)\alpha(\xi){\partial\alpha\over\partial T}(\xi)\eta\in L^{p/(p-2) }(\Omega)$ is continuous
thus the second term in $I_{1,m}$ passes to the limit as $m$ tends to infinity.

Analogously, we take to the limit in $I_{2,m}$
 observing that the strong-weak convergence product
$\nabla{\xi}_m\cdot\nabla\phi_m\rightharpoonup \nabla\xi\cdot\nabla\phi$ holds
 in ${ L}^{p/2}(\Omega)$.
 
In order to be in conditions for finding that $\theta$ is a limit solution,
let us prove
the continuity of the solution mapping $\xi\in
W^{1,p}(\Omega)\mapsto\phi=\phi(\xi)\in  W^{1,s}(\Omega)$  in
the strong topology  for 
$s=2pn/(p+n)<p$. 
Take the difference of (\ref{auxphi})$_m$ and (\ref{auxphi})
verified by the solutions $\phi_m$ and $\phi$, respectively, and
choose $\eta=\phi_m-\phi$ as a test function. Thus,
it results
\begin{eqnarray*}
\sigma_\#\|\nabla (\phi_m-\phi)\|_{2,\Omega}^2\leq 
\int_\Omega (\sigma(\xi)-\sigma(\xi_m))\nabla\phi\cdot\nabla(\phi_m-\phi) dx+\\
+\int_\Omega (\sigma(\xi)\alpha(\xi)\nabla\xi-\sigma(\xi_m)\alpha(\xi_m)\nabla\xi_m)
\cdot\nabla(\phi_m-\phi) dx\longrightarrow 0,
\quad\mbox{as $m\rightarrow\infty$.}
\end{eqnarray*}
Then, we conclude that $\nabla{\phi}_m \rightarrow\nabla\phi$ in ${\bf L}^{2}(\Omega)$,
and consequently $\nabla{\phi}_m \rightarrow\nabla\phi$ a.e. in $\Omega$ and
$|\nabla{\phi}_m|^2 \rightarrow |\nabla\phi|^2$ in ${L}^{s/2}(\Omega)\hookrightarrow L^{pn/(p+n)}(\Omega)$.
Thus $I_{3,m}$ passes to the limit as $m$ tends to infinity,
 concluding the proof of weak continuity of the operator $\mathcal L$.
 
Then the Schauder fixed point theorem can be used and it guarantees the
existence of $(\theta,\phi)$ in the conditions to Theorem
\ref{exist}. 

\subsection{The validation of the estimates (\ref{cotaphi})-(\ref{cotau})}

Let  $(\theta,\phi)\in W^{2,pn/(p+n)}(\Omega)\times W^{1,p}(\Omega)$ be a weak solution
to the variational problem (\ref{heatj})-(\ref{pbphij}).

Choose $\eta=\phi\in W^{1,p}(\Omega)$ as a test function in (\ref{pbphij}).
Using (\ref{smin}), the upper bound of $|\alpha|$ and the Sobolev-Poincar\'e inequality 
 then (\ref{cotaphi}) holds.
 
From the regularity theory for the Robin-Laplace problem
and by virtue of the existence of a solution $\theta\in W^{2,pn/(p+n)}(\Omega)$ 
we proceed as in (\ref{cotatp}) now for $k>\alpha_\#h_\#$
resulting the estimate
\[
\| \nabla\theta\|_{p,\Omega}
\leq K(A\|\nabla\theta\|_{p,\Omega}^2+B\|\nabla\theta\|_{p,\Omega}
\|\nabla\phi\|_{p,\Omega}+\sigma^\#\|\nabla\phi\|_{p,\Omega}^2+
\|g\|_{p/2,\Omega}),
\]
 with $K$ denoting a constant independent on $k$.
Combining this result with the estimate (\ref{cotap}) with $\xi$ replaced by $\theta$ and using (\ref{ss}) 
we conclude (\ref{cotaup}).

Choose $\eta=\theta\in W^{1,p}(\Omega)$ as a test function in (\ref{heatj}).
Then applying the H\"older inequality and using the assumptions (H1)-(H4) it follows
\begin{eqnarray*}
k\|\nabla \theta\|^2_{2,\Omega}\leq\alpha^\#h_1
\|\theta\|^{2}_{2,\Sigma}+
\Big(A\|\nabla\theta\|_{2,\Omega}^2+B\|\nabla
\theta\|_{2,\Omega}\|\nabla\phi\|_{2,\Omega}+\\+\sigma^\#\|\nabla
\phi\|_{2,\Omega}^2+
\|g\|_{1,\Omega}\Big)\|\theta\|_{\infty,\Omega}.
\end{eqnarray*}

This yields the estimate (\ref{cotau}).

\section{Proof of Theorem \ref{exist1}}
\label{sexist1}

 For each given $k>0,$ let $(\theta_k,\phi_k)$ be a solution to
(\ref{heatj})-(\ref{pbphij}) in accordance with Theorem \ref{exist}. From
estimates (\ref{cotau}) and (\ref{cotaphi}) there exist subsequences
still denoted by $\theta_k$ and $\phi_k$ such that, for
$k\rightarrow+\infty,$
\begin{eqnarray*}
\nabla \theta_k\rightarrow 0 
\quad\mbox{ in }{\bf L}^2(\Omega);\qquad
\theta_k\rightharpoonup \Theta \quad\mbox{ in }
W^{1,p}(\Omega)\hookrightarrow\hookrightarrow C(\bar\Omega);\\
\phi_k\rightharpoonup\phi\qquad \mbox{ in $H^1(\Omega)$},
\end{eqnarray*}
with $\Theta$ constant on $\bar\Omega$.
Hence we can pass to the limit  in (\ref{pbphij}) as $k$ tends to infinity resulting
\[
\int_\Omega\sigma(\cdot,\Theta)\nabla\phi\cdot \nabla\eta dx=
\int_{\partial\Omega}
h\eta ds,\quad\forall \eta\in W^{1,p'}(\Omega),
\]
or equivalently (\ref{t1})-(\ref{t2}).

In particular, if we take $\eta=1$  in (\ref{heatj})
we can pass to the limit as $k$ tends to infinity resulting (\ref{tt}).

Using (\ref{smin}) it follows
\begin{eqnarray*}
\sigma_\#\int_\Omega|\nabla\phi|^2dx+\int_\Omega gdx
\leq \Theta\int_{\partial\Omega}\alpha(\cdot,\Theta)h ds\leq
\sigma^\#\int_\Omega|\nabla\phi|^2dx+\int_\Omega gdx.
\end{eqnarray*}
Taking into account that the assumption (H4)$_+$ or (H4)$_-$ and also (\ref{c1})
imply
\begin{eqnarray*}
\int_{\partial\Omega}\alpha(\cdot,\Theta)h ds &\leq &\alpha^\# 
\int_{\partial\Omega}|h|ds;\\
\int_{\partial\Omega}\alpha(\cdot,\Theta)h ds &\geq &
\alpha_\#h_\#|\Gamma|+\alpha^\# h_1|\Sigma|>0,
\end{eqnarray*}
we derive (\ref{tt1}), concluding the proof of Theorem \ref{exist1}.

\section{The two-dimensional limit model}
\label{2d}

\subsection{Proof of Theorem \ref{teo1}}

Arguing as in Theorem \ref{exist},
the Schauder fixed point argument can be applied.
For this, it is sufficient to see that
the  regularity relative to
 $W^{1,p}_{\rm loc}(\bar\Omega)=W^{1,p}(\Omega)$
($\Omega$ bounded)
 can be applied for the unique   weak solution of the variational problem (\ref{auxphi})
 for $p=2+\epsilon>2=n$ in accordance to Proposition \ref{pphi}.
 Thus  Proposition \ref{propt2}
 guarantees the existence of $\theta\in W^{2, 2p/(p+2)}(\Omega)$
 verifying (\ref{cotatp}).
For every $\xi\in \bar B_R$, inserting (\ref{cotaphie}) into (\ref{cotatp}) it follows
\begin{eqnarray*}
\| \theta\|_{2,2p/(p+2),\Omega}
\leq K(AR^2+BR {K_2}(\sigma^\#\alpha^\#R+\|h\|_{2,\partial\Omega} ) 
+\\
+\sigma^\# {K_2}^2(\sigma^\#\alpha^\#R+\|h\|_{2,\partial\Omega} ) ^2+
\|g\|_{p/2,\Omega}).
\end{eqnarray*}
  Next arguing as in Section \ref{sexist}
it leads to a smallness condition.
 From the continuity of the Nemytskii operator $\sigma $ due to  the Krasnoselski Theorem
we can proceed as in Section \ref{swc},
considering that $I_{3,m}$ passes to the limit as $m$ tends to infinity since
\begin{eqnarray*}
\sigma(\xi_m)\rightharpoonup \sigma(\xi)\qquad \mbox{weakly* in }L^\infty(\Omega);\\
|\nabla{\phi}_m|^2\eta \rightarrow |\nabla\phi|^2\eta\qquad \mbox{in }L^{1}(\Omega).
\end{eqnarray*}
This concludes that we are in the conditions of applying Theorem \ref{fpt}.

\subsection{Proof of Proposition \ref{exist2}}

We can proceed as in Section \ref{sexist1} considering the existence of the sequence of solutions
is provided by Theorem \ref{teo1} and the proposition follows.

\end{document}